\documentclass[reqno,a4paper,11pt]{amsart}

\usepackage{color}
\usepackage[latin1]{inputenc}
\usepackage[english]{babel}
\usepackage{amsmath,amssymb,amsfonts,amsthm,enumerate}
\usepackage[T1]{fontenc}
\usepackage[active]{srcltx}
\usepackage{epsfig}
\usepackage{graphicx}

\numberwithin{equation}{section} 

\newcommand{\tauV}{{\kern-3pt\tau}}

\newcommand{\oVVVk}{\overline{\mbox{\boldmath$V$}}\kern-3pt}

\newcommand{\tVVVk}{\tilde{\mbox{\boldmath$V$}}\kern-3pt}

\begin{document}
\title[Free Boundary Problems]
{Free Boundary Problems:\\[2mm]
The Forefront of Current and Future Developments}

\author[Chen]{Gui-Qiang Chen}
\address{Mathematical Institute, University of Oxford, Oxford, OX2 6GG, UK}
\email{chengq@maths.ox.ac.uk}

\author[Shahgholian]{Henrik Shahgholian}
\address{Department of Mathematics, KTH Royal Institute of Technology, 100~44 Stockholm, Sweden}
\email{henriksh@kth.se}

\author[Vazquez]{Juan-Luis Vazquez}
\address{Department of Mathematics, Universidad Aut\'{o}noma de Madrid,
28049 Madrid, Spain}
\email{juanluis.vazquez@uam.es}
\maketitle

\section{Introduction}

\medskip
The term  {\it Free Boundary Problem {\rm (FBP)}} refers, in the modern
applied mathematical literature, to a problem in which one or several variables
must be determined in different domains of the space, or space-time,
for which each variable is governed in its domain by a set of state laws.
If the domains were known, the problem reduces to solving the equations,
usually ordinary differential equations (ODEs) or partial differential equations (PDEs).
Now, the novelty of FBPs lies in the fact that the domains are {\it a priori} unknown
and have to be determined as a part of the problem, thanks to a number of
free boundary conditions that are derived from
certain physical laws or other constraints governing the phase transition.

\smallskip
The interplay of the state laws for the single phases
and the special phase transition conditions leads to
a mathematical discipline that combines analysis and geometry in sophisticated ways,
together with mathematical modeling based on physics and engineering.
Let us point out a basic division of the set of free boundary problems:
In the processes that evolve in time, one of the main difficulties of the theory is to
track the movement of the free boundary or boundaries; these problems are also called
{\sl moving free boundary problems}.
The processes where time does not appear have the common name, {\sl FBP.}
Another basic distinction is between several-phase problems and
the simpler one-phase problems,
where the free boundary is the bounding hypersurface of the phase that is governed
by the equations. Normally, a one-phase problem can be considered as a two-phase problem
with a trivial second phase, usually the vacuum.

\smallskip
It is probably very hard, if not impossible, to track the origins of free boundary problems.
The origin of the modern discipline owes much the famous {\it Stefan problem} that describes
the joint evolution of a liquid and a solid phase, a question considered in 1831
by G. Lam\'{e} and B. Clapeyron in relation to the problems of ice formation in the polar seas.
The problem is named after J. Stefan \cite{Stefan},
who introduced the general class of
such problems around 1890, and performed both numerical and theoretical studies of
the phase-transition problems related to heat transfer.
The rigorous analytical treatment exceeded the tools of the time, and the difficulty may be
understood from the fact that the basic existence and uniqueness result was only proved
by S. Kamenomostkaya \cite{Kam} and O. Oleinik \cite{Oleinik},
once the notion of weak solutions and nonlinear analysis were ready.

\smallskip
In the twentieth century, new directions were developed, and problems concerning
free surfaces flows, shock waves, and water waves became of central importance.
During early 50's until early 70's, several mathematicians
(such as P. Lax, H. Lewy, J.~L. Lions, G. Stampacchia, among others)
contributed to the growth of the field, and eventually the weak and variational formulation
and techniques for these problems were fully developed.
The {\it obstacle problem} became a basic study case in the mathematical theory.
There are many names here unmentioned, but not forgotten.
The book by G. Duvaut and J.~L. Lions was probably one of the first to make an extensive list
of problems available to mathematical community.
At that stage of the development of the topic, the mathematical treatment was anchored in
the new tools of nonlinear functional analysis and PDEs,
and several deep contributions were made to obtain the regularity of solutions
and to develop the theory towards other directions.

The 70's and 80's saw a continuous progress in the treatment of an increasing number of models,
some of them became classical in the analysis of FBPs:
the dam problem, the plasma problem, the porous medium equation,
the $p$-Laplacian equation, the Hele-Shaw problem, the Muskat problem, among many others.
At the same time, the difficult task of understanding the finer points of the geometry
of the free boundaries become central. The work of L. Caffarelli allowed to make
a remarkable progress in the questions of regularity of FBPs after his seminal paper \cite{Caffarelli77}.
It is also the time when a mathematical FBP community was formed, a regular series of international FBP meetings
was scheduled, and all this served to maintain the impulse and activate the innovation.
With time, different books  (\cite{Fried1982, Elliot-book, Meirmanov-book, CSbook, vazquezPME, PSU} to  mention a few)
were written to consolidate the progresses that had been done in different areas.
See also the survey paper  \cite{Fried2000} reporting on the situation around year 2000.

\smallskip
Recent decades have witnessed a rapid widening of the subject area by incorporation
of important free boundary topics coming from different areas:
In Finance, FBPs appear to determine the optimal  exercise value  in Black-Scholes models;
in Mathematical Biology, they indicate the moving fronts of populations or tumors;
many new and emerging FBPs have arisen from Fluid Mechanics; free boundaries appear in aggregation/swarming processes;
in Geometry, much attention is now given to geometrical flows,
including the curvature flows where different types of curvatures involved.


\medskip

\section{Current and Future Developments}

Today, the study of FBPs is intensely pursued from various aspects
(experimental, numerical, and theoretical),
the subject is continuously finding new grounds for applications,
and new fundamental theoretical questions continue to emerge.
These developments, in particular, ask for new analytical and numerical methods,
and improvements of existing algorithms and tools
to handle extremely complex problems.

\smallskip
This theme issue contains several articles that represent the mainstream
in today's and future directions of FBPs (even though not exhaustive).

\subsection{Theoretical developments}

From the current status of theoretical directions, we have selected a number of topics that will be covered in this issue.

\smallskip
\noindent
{\sl Nonlocal Phenomena:}
Many recent developments involve phenomena for which the diffusion process is of a nonlocal
nature. Thus, in continuum mechanics and fluid dynamics, owing to the presence of many scales such as in polymers,
there exists a global interaction through force fields,
or induced by spatial phenomena in cases of surface diffusivity.
Indeed, if heat flow is blocked by an insulated wall,
heat will propagate along the wall by the influence of the internal heat,
through the Dirichlet to Neumann nonlocal diffusion kernel.
This is the case of the quasi-geostrophic equation or models of planar
crack propagation in viscoelastic solids.
Nonlocal diffusion behaviour also takes place in probability, when the random kicks affecting the particles
cease to be {\it infinitesimal} and become {\it jump processes} (L\'evy processes).
Many phenomena involving phase transitions and free boundaries re-appear in this context.
Nonlocal behaviour has many aspects and scales (thin membrane, thick membrane, different scales
for interior and boundary diffusivity, {\it etc.}) and presents many new mathematical difficulties.
It plays a crucial role in applications such as elasticity (the Signorini problem),
optimal insulation, and mathematical finance (American option with L\'evy process).

It is also noteworthy that,  in the last decade, there has been a very definite progress in the mathematical
understanding of nonlocal problems of stationary or evolution type, with or without free boundaries.
Some important papers in the latter  direction are \cite{CaffSilv2007, CaffVass2010, CaffVaz2011, Vaz2012abel}.
More references will be mentioned below.

\smallskip
\noindent
{\sl Two- and Multi-Phase Problems:}
Many phenomena modeled by PDEs give rise to several phases of behaviour,
in which two- or multi-phases of free boundaries are natural.
They may arise as limiting cases of solutions to systems of equations.
Several phases in flow problems also appear in applied problems such as in oil related
industry (where several phases of materials meet or mix), {\it e.g.}, oil, gas, air, saltwater,
sand, {\it etc.} The well-known Muskat problem is an example of such phenomena, see \cite{CCFGL}.
Another problem is the so-called segregation in reaction-diffusion problems
which models a competition between several species \cite{Mat-Mim,DHMP}.
The theory of these problems are still at an embryonic stage of development; see \cite{CTV}
 for latest developments concerning segregation problems.
We refer also to \cite{PSU} and the references therein, for other types of two-phase problems related to
 obstacle problem.

\smallskip
\noindent
{\sl Combinatorial Aspects and Random Walks:}
Many applied problems result in particle interactions in various forms.
Such problems, driven by exterior or interior forces, chemical reactions,
or population growth, have been studied to some extent within applied sciences.
When particles move ({\it e.g.} cell growth, crystallisation, discrete flows) randomly
and occupy regions, they interact with other particles in their environment.
Such interaction can be of various forms:
They can be annihilated, or annihilate other particles;
they can also freeze or vaporise and behave in many other ways.
Such movements usually are described by random walks,
and the density function representing the population, or the amount of chemical substances,
are given by discrete harmonic functions.

Models for one-phase version of these discrete problems have
been developed
by several mathematicians. We refer to one of the latest results, \cite{JLS}, and the references cited therein.
The problem has some history in combinatorial problems
such as chip firing, random algebraic sums (Smash-sum), and colouring problems.
The two- and multi-phase versions of such problems have yet not been developed,
even though there are huge amount of corresponding models.

\smallskip
\noindent
{\sl Shock Waves:}
\noindent
The study of solutions with shocks in various models of compressible fluid dynamics
leads to well-known free boundary problems.
One class of such problems involves steady and self-similar solutions of
the multidimensional Euler equations.
Shocks correspond to one kind of discontinuities in the fluid velocity, density, and  pressure,
which are discontinuities mathematically in the solutions or in their derivatives
depending on the model equations ({\it cf.} \cite{Dafermos-book,CFeldman-A}).
Important physical applications involving shocks include self-similar shock reflection-diffraction
and steady shocks in nozzles and around wedges/cones or airfoils.
Such shocks are important in the mathematical theory of multidimensional conservation laws
since steady/self-similar solutions with shocks are building blocks and asymptotic attractors
of general solutions ({\it cf.} \cite{CFeldman2,CFeldman-A}).
In recent years, mathematical progress has been achieved in several long-standing problems,
especially shock reflection-diffraction problems,
for potential flow in which many different flow regimes are possible.
Future challenges concern:

\smallskip
$\bullet$ Complicated patterns of regular shock reflection-diffraction and other configurations for potential flow;

$\bullet$   Extension of the results of shock reflection-diffraction, including regular/Mach configurations,
    to the full Euler equations;

$\bullet$  Construction of steady solutions in a nozzle ({\it e.g.} the de Laval nozzle),
or around a non-straight wedge/cone or an airfoil;

$\bullet$  Analysis of free boundaries in solutions to the multidimensional Riemann problem.

\smallskip
Other related discontinuities (as free boundaries) in compressible fluid mechanics
include vortex sheets and entropy waves ({\it cf.} \cite{CFeldman-A,CWang1}).


\subsection{Numerical developments}

In the early days, numerical approaches to FBPs were frequently ad hoc based,
say, on trial free boundary methods or coordinate transformations.
Systematic approaches based on variational methods, variational inequalities, \cite{Kind-Stam-book}
and weak formulations of degenerate nonlinear parabolic equations were developed
as the mathematical theory became available. In recent years, important
developments include the following:

\smallskip
$\bullet$ Diffuse interface models (such as the phase field and Cahn-Hilliard models)
  with applications to curvature flows, solidification and phase transformations
  in material science;

$\bullet$ Level set methods for evolving fronts including applications to fluid flow and
     image processing;

$\bullet$ Variational front tracking methods for geometric PDEs; for instance, interfaces
      involving curvature effects (such as surface tension and bending);

$\bullet$ Extensive mathematical contributions to the stability, well-posedness and
     rigorous error analysis of discrete approximations to free boundary problems
      and degenerate nonlinear elliptic and parabolic equations.

$\bullet$ Adaptive methods appropriate for free boundary and interface problems.

\smallskip
The need to simulate ever-more complicated large systems, stemming from the increase
in computing power and the development of computational tools, stimulates new
questions and problems for analysts. Therefore, there is a natural symbiosis
between analysts, modelers, and computational mathematicians.
Increased computing power together with the demand of applications
have led to the study of systems of PDEs in domains with complex morphology.
These complex multi-physics models frequently involve interfaces and free boundaries.
Although one may obtain detailed information about sub-problems (say, the obstacle problem)
or local-in-time existence results by using analytical techniques,
the full complex system is required to be simulated in scientific and engineering applications.
Nonlinear degenerate PDEs with interfaces and free boundaries are
notoriously difficult to solve numerically. The
best
numerical
methods depend on good analytical approaches, and sometimes they
promote new advances in the PDE theory, such as the formulation of mean
curvature flows beyond the onset of singularities. Future trends in
numerics may include the following:

\smallskip
$\bullet$
Surface finite elements for geometric PDEs, surface processes on
interfaces;

$\bullet$ Diffuse interface and phase field methods for two phase flow with
surfactants, phase transformations in materials;

$\bullet$ Numerical optimization of free boundaries for control and inverse
problems;

$\bullet$ Burgeoning applications in biology and medicine;

$\bullet$ Numerical methods for fully nonlinear equations;

$\bullet$ Numerical methods in homogenization, random media and random
surfaces including stochastic equations;

$\bullet$ Computational methods for free boundary problems for shock waves,
vortex sheets, entropy waves, and related compressible flows;

$\bullet$ Adaptivity (mesh refinement, coarsening and smoothing) for
surfaces including topological change.

\subsection{Applications and the present volume}

As indicated above, the study of FBPs is an extremely broad field
due to the abundance of applications in various sciences and
real world problems, including physics, chemistry, engineering,
industry, finance, biology, and other areas.
This theme issue is also a showcase of applications of FBPs
with a wide coverage of different areas from fluid mechanics
to biology/medicine to finance.

\smallskip
In \cite{Elliott},  A. Alphonse and C. Elliott formulate a Stefan problem on an evolving
hypersurface and study the well-posedness of weak solutions given $L^1$ data.
As mentioned earlier, the Stefan problem is the prototypical time-dependent free boundary problem.
It arises in various
forms in many models in the physical and biological sciences. In order to do their analysis,
they have to first develop function spaces and results to handle the equations on evolving surfaces
in order to give a natural treatment of the problem. The applications to numerical treatment of applied
problems are natural and interesting.

D. Apushkinskaya and N. Uraltseva in \cite{ApushkinskayaUraltseva-A}
study the regularity of free boundaries in problems with hysteresis that arise in  modeling  biological
and chemical processes "with memory". Such models lead to two-phase free boundary problems.

In \cite{Bucur}, D. Bucur and Velichkov show
how the analysis of a general shape optimization
problem of spectral type can be reduced to the
analysis of particular free boundary problems. The paper presents an overview of recent developments.

L. Caffarelli and H. Shahgholian in \cite{caff-sh} give a survey on the approaches to regularity
theory for four different free boundary problems,
which are  obstacle, thin obstacle, minimal surfaces and cavitation problems.

The paper by J. A. Carrillo \& J. L. V\'azquez  \cite{CarrilloVazquez-A}
is devoted to problems of diffusion or aggregation. Of concern is the interaction
between nonlocal diffusion with long-range interactions and nonlinearities of the degenerate
type that may slow the movement. The question is then whether there is a free boundary.
The answer depends on the models in a delicate way.
Once this settled, the regularity and asymptotic behaviour take the scene.

The paper by G.-Q. Chen \& M. Feldman \cite{CFeldman-A} provides a survey
of recent activities in deriving and analyzing free boundary problems
in shock reflection/diffraction and related transonic flow problems.
Shock waves are steep fronts that propagate in the compressible fluids
which are fundamental in nature, especially in high-speed fluid flows.
It is shown that several longstanding shock reflection/diffraction
problems can be formulated as free boundary problems,
some recent progress in developing mathematical ideas, approaches,
and techniques for solving these problems is discussed.

P. Constantin investigates in \cite{Constantin-A}  the dynamics of vortex patches
in the Yudovich phase space. He derives an approximation for the evolution of the vorticity
in the case of nested vortex patches with distant boundaries, and studies its long time behaviour.

In \cite{Cordoba}, D. C\'ordoba, J. G\'omez-Serrano, and A. Zlatos
deal with a famous problem in the theory of flows in porous media, the Muskat Problem.
They show that there exist solutions of the  problem that shift stability regimes:
They start unstable, then become stable, and finally return to the unstable regime.

A. Figalli and H. Shahgholain in \cite{Fig-Sh} present a survey on unconstrained free boundary problems involving
elliptic operators. The main objective is to discuss a unifying
approach to the optimal regularity of solutions to
these  matching problems.

Biological applications have attracted enormous attention lately,
because of their importance to society. A. Friedman in \cite{Friedman-A} reviews several free boundary
problems that arise in the mathematical modeling of biological processes. The biological topics are
quite diverse: cancer, wound healing, biofilms, granulomas, and atherosclerosis. For each of these
topics there is a description of  the biological background and the mathematical model.

J. Glimm et al. in \cite{GLKH} review the existence and non-uniqueness
for the Euler equations of fluid flow. The non-uniqueness of solutions is a fundamental
obstacle to scientific predictions, preventing
the success of such goals as predictive science, validation of solutions and quantification of the
uncertainties on which engineers base their design conclusions. Some of the non-uniqueness
of the Euler equations appears to be mathematical in nature, and it may be assumed to disappear after some future,
deeper level of mathematical results have been obtained.
However, some of the nonuniqueness is physical, or at least numerical in origin,
due to the underspecification of the Euler equations on physical grounds.
A mitigating strategy is proposed in the paper.

M. Hadzic and S. Shkoller in \cite{Shkoller} study the stability of steady states in the Stefan problem
for general boundary shapes, which is a contribution to the theory of one of the most classical examples of
the free boundary field. They prove the global-in-time stability of the steady states of the classical one-phase Stefan
problem without surface tension, assuming a sufficient degree of smoothness on the initial domain.

R. Nochetto and  R. Otarola  in \cite{Nochetto}  review the finite element approximation of the
classical obstacle problem in energy and max-norms and derive error estimates for both the solution and
the free boundary. They  present an optimal error analysis for the thin obstacle problem.
Finally, the authors discuss the localization of the obstacle problem for the fractional Laplacian
and prove quasi-optimal convergence rates.

In the same area, B. Perthame and N. Vauchelet in \cite{Perthame} deal with the incompressible limit of a
mechanical model of tumor growth with viscosity.  In particular, various mathematical models of tumor growth are available in the
literature.  A first class of models describes the evolution of the cell number density when considered as a continuous
visco-elastic material with growth. A second class of models describes the tumor, as a set and rules for the free
boundary are given related to the classical Hele-Shaw model of fluid dynamics.


\bigskip
\bigskip
{\bf Acknowledgments:}
G.-Q. Chen's research was supported in
part by  the Royal Society--Wolfson Research Merit Award (UK)
and the UK EPSRC Award to the EPSRC Centre for Doctoral Training
in PDEs (EP/L015811/1).
H. Shahgholian's research was supported in part by  Swedish Research Council.
J.-L. Vazquez's research was supported in part by  the Spanish Research Project MTM2011-24696.
The authors would like to thank the Isaac Newton Institute for Mathematical Sciences,
Cambridge, for support and hospitality during the 2014 Programme on
{\it Free Boundary Problems and Related Topics} where this work  was undertaken.

\bigskip


\end{document}